%% file: ConLaw_Paper_1_ArXiV.tex
\documentclass[reqno,10pt]{amsart}
\usepackage{empheq}
\numberwithin{equation}{section}
\input{./common/header_arxiv.tex}

\usepackage{pgfplots}
\pgfplotsset{compat=1.17}
\usepgfplotslibrary{fillbetween}
\usepackage{float}
\allowdisplaybreaks

    \title{A different derivation of conservation laws for water waves}
    \author{Katie Oliveras$^\dagger$ and Salvatore Calatola-Young$^\ddagger$}
    \date{\today}
    \address{$^\dagger$ $^\ddagger$Mathematics Department, Seattle University, Seattle, WA 98122}
    \email{$^\dagger$oliveras@seattleu.edu}
    \urladdr{$^\dagger$fac-staff.seattleu.edu/oliveras/web} 
\begin{document}
    \begin{abstract}
        We consider a new nonlocal formulation of the water-wave problem for a free surface with an irrotational flow based on the work of Ablowitz, Fokas, and Musslimani \cite{ablowitz2006new} and presented in \cite{NonlocalNonlocal}. The main focus of the short paper is to show how one can systematically derive Olver's eight conservation laws not only for an irrotational fluid, but also for constant vorticity (linear shear flow) without explicitly relying on the underlying Lie symmetries. This allows us to make draw new conclusions about conservation laws and posit the existence of additional, nonlocal, conservation laws for the water-wave problem.
    \end{abstract}
\maketitle

\section{Introduction}
    \input{./sections/introduction.tex}
\section{Prior Work Related to Conserved Densities} \label{sec:background}
    \input{./sections/overview.tex}
\section{The Nonlocal/Nonlocal Formulation} \label{sec:equations}

\input{./sections/equations.tex}
\section{Conserved Densities for an Irrotational Fluid}\label{sec:conlaw_irrotational}
    \input{./sections/conservation_laws_1D_surface.tex}
\section{Conserved Densities for Linear Shear} \label{sec:conlaw_conVort}
    \input{./sections/conservation_laws_linear_shear.tex}    
\section{Acknowledgements}
    \input{sections/acknolwledgements.tex}
\input{sections/appendix.tex}
\bibliographystyle{alpha}
\bibliography{./common/references.bib}
\end{document}

%% file: common/header_arxiv.tex
\usepackage[margin=.75in]{geometry}
\usepackage[usenames,dvipsnames]{xcolor}
\definecolor{linkColor}{HTML}{f10000}
\usepackage[    colorlinks  =   true,
                linkcolor   =   Red,
                urlcolor    =   NavyBlue,
                citecolor   =   NavyBlue,
                anchorcolor =   Emerald]{hyperref} 

\usepackage{amsmath, amsthm,mathrsfs}
\usepackage{booktabs}
\usepackage{array}
\usepackage{cleveref}
\usepackage{setspace}
\usepackage{soul}
\usepackage{cite}
\usepackage{graphicx}
\newtheoremstyle{indented}{3pt}{3pt}{\addtolength{\leftskip}{1.5em}
        \itshape}{}{\bfseries}{.}{.5em}{}

\theoremstyle{indented}
\newtheorem{remark}{Remark}
\AtBeginDocument{
  \addtocontents{toc}{\small}
  \addtocontents{lof}{\small}
}
\newcommand{\dt}{\frac{d}{dt}}

\newcommand{\dx}{\frac{d}{dx}}

\newcommand{\grad}{\nabla}
\newcommand{\D}{\mathscr{D}}
\newcommand{\vn}{\vec{n}}
\newcommand{\vt}{\vec{t}}
\newcommand{\dotn}{\cdot\vn}
\newcommand{\dott}{\cdot\vt}

\newcommand{\cint}[1]{\oint_{\partial\D}\left[#1\right]\,ds}
\newcommand{\sint}[1]{\int_{\mathscr{S}}{\left[#1\right]\,dx}}
\newcommand{\bint}[1]{\int_{\mathscr{B}}{\left[#1\right]\,dx}}
\newcommand{\wint}[1]{\int_{-\infty}^\infty{\left[#1\right]\,dx}}

\newcommand{\dint}[1]{\iint_{\mathscr{D}}#1\,dx\,dz}

\newcommand{\odint}[1]{\oint_{\partial\mathscr{D}}\left(#1\right)}

\usepackage{fontawesome}
\usepackage{tcolorbox}
\tcbset{width=\textwidth, boxrule=1pt}

%% file: sections/introduction.tex
In \cite{olver}, Benjamin and Olver derived the symmetries corresponding to the free-boundary problem for an inviscid and irrotational fluid for both a one- and two-dimensional free surface. To find these conserved densities, the authors use transformation groups and prolongation methods to recover the symmetry group.  From these symmetries, the conserved densities are derived.  

For a one-dimensional free surface without surface tension, eight conserved densities were found.  In a subsequent paper, Olver proved that these eight conservation laws (or seven in the presence of surface tension) are the only conservation laws that arise from the infinitesimal symmetries for the free boundary problem.  To the best of our knowledge, the work of Benjamin and Olver represent the only known systematic derivation of the appropriate conserved densities.   

In this short note, we describe another systematic derivation of these conserved densities via a weak formulation of the problem.  In \Cref{sec:background}, we present the equations of motion as well as the notion of a conserved density for the water-wave problem. We then consider the nonlocal/nonlocal formulation presented in \cite{NonlocalNonlocal} for an inviscid, irrotational fluid posed on the whole line with a one-dimensional free surface.  After a brief introduction to the formulation in \Cref{sec:equations}, we illustrate a process for deriving Benjamin \& Olver's conservation laws without \emph{explicitly} using the symmetries.  

In \Cref{sec:conlaw_irrotational}, we reconstruct all of Benjamin \& Olver's conserved densities with a special note about the Hamiltonian. In \Cref{sec:conlaw_conVort}, we extend the results to constant vorticity showing that Olver's conservation laws are easily extended to linear-shear - though the results conserved densities in terms of surface variables are nonlocal - to the best of our knowledge, these results are new.  

While Olver's proof \cite{olver1983conservation} establishes that there are exactly eight non-trivial conservation laws for a one-dimensional surface, to the best of our knowledge, the question of additional nonlocal conservation laws remains unknown.  

%% file: sections/overview.tex
\begin{figure}[htp]
        \centering
        \begin{tikzpicture}
            \begin{axis}[
                            hide axis,
                            xmin        =   -13,
                            xmax        =   13,
                            ymin        =   -1.15,
                            ymax        =   .8,
                            domain      =   -10:10,
                            width       =   \textwidth,
                            height      =   1.5in,
                            samples     =   300,
                            black,
                            thick
                        ]
                \addplot[name path = wave] {.35*2^(-.2*x^2)} node [pos=1, right] {\(z = \eta(x,t)\)}; 
                \addplot[name path = bottom]{-1} node [pos = 1, right] {\(z = -h\)};
                \addplot[name path = floor] {-1.15}; 
                
                \addplot[fill = NavyBlue!30] fill between[of=wave and bottom];
                \addplot[fill=Brown!60] fill between[of=floor and bottom];

            \end{axis}
        \end{tikzpicture}
        \caption{Fluid Domain}\label{fig:fluidDomain}\end{figure}
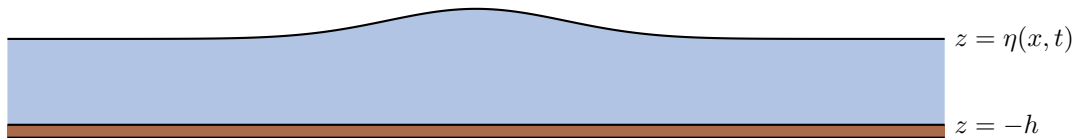

    We consider the following equations for the motion of an inviscid, irrotational fluid problem posed on the whole line as illustrated in \Cref{fig:fluidDomain}:
        \begin{align}
            &\phi_{xx} + \phi_{zz} =0,  & &(x,z)\in\mathscr{D}, \label{eqn:laplace1d} \\ 
            &\phi_t + \frac{1}{2}\vert\nabla\phi\vert^2 + gz + p/\rho = 0,  & &(x,z)\in\mathscr{D},\label{eqn:bernoulliBulk} \\  
            &\phi_z =0, & &z = -h, \label{eqn:kinematicBottom1d}  \\
            &\eta_t + \phi_x\eta_x =\phi_z,  &&z = \eta(x,t),\label{eqn:kinematic1d}  \\
            &\phi_t + \frac{1}{2}\vert\nabla\phi\vert^2 + g\eta =0, && z = \eta(x,t), \label{eqn:dynamic1d} 
        \end{align}
    where \(\phi\) is the velocity potential, \(\eta\) is the free-surface, and \(p\) is the fluid pressure. The physical parameters \(g\) and \(\rho\) represent gravity and the constant fluid density respectively.  
 
    In what follows, we consider a simpler scenario that what was prescribed in \cite{olver}.  Specifically, we consider a fluid domain as shown in \Cref{fig:fluidDomain} where we define the following: 
    \begin{displaymath}
        \mathscr{D} = \lbrace (x,z):x\in\mathbb{R}, -h<z<\eta(x,t)\rbrace, \qquad\mathscr{S} = \lbrace (x,z):\mathbb{R}, z = \eta(x,t)\rbrace,\qquad \mathscr{B} = \lbrace (x,z):x\in\mathbb{R}, z = h\rbrace,
    \end{displaymath}
    with the path \(\Gamma\) defined such that \(\partial\mathscr{D} = \Gamma \cup\mathscr{S}\).  

    In this paper, we present results for waves posed on a domain where the free-surface \(z = \eta(x,t)\) is single-valued, continuously differentiable, and decays sufficiently fast as \(\vert x \vert \to \infty\). Furthermore, we take the mathematicians view-point and impose that \(\phi(x,z,t)\to 0\) sufficiently fast as \(\vert x \vert \to \infty\). The problem is easily relaxed to a variety of configurations including parametric representations for the free-surface and by relaxing the restriction \(\phi \to 0\) to become \(\left\vert\grad\phi\right\vert \to 0\).  While there are some minor changes in the formulation, the results presented here persist. 

    \subsection{Conserved Density Forms}
        
        Before we proceed to generate the conservation laws, it is important to note precisely what is meant by conservation law for this free-boundary problem.  Following the work of \cite{olver,olver1983conservation}, a function \(T\) depending on \(x, t, \eta, u,\) and the derivatives of \(\eta\) and \(u\) over the free-surface \(\mathscr{S}\) is a \emph{conserved density} for the free-boundary problem if there exists a vector function,
        \[
            \vec{\mathbf{F}} = \begin{bmatrix}F_1\\F_2\end{bmatrix},
        \] 
        with \(F_1\) and \(F_2\) depending on \(x, z, t, u,\) and the derivatives of \(u\) in the region \(\mathscr{D}\), and a function \(W\) depending on \(x, t, \eta, u,\) and their derivatives on \(\mathscr{S}\) such that for all solutions \(u = f(x,z,t),\eta = g(x,t)\) of the free boundary problem, 
        \[
            \mathbf{D}_t T = \vec{\mathbf{F}}\dotn + \mathbf{D}_x W \text{ on } \mathscr{S},
        \] with 
        \[
            \grad\cdot\vec{\mathbf{F}} = 0\quad \text{ in }\quad  \mathscr{D}.
        \]

        As described in \cite{olver1983conservation}, this implies that conservation laws for the free-surface \(z = \eta\) take the general form 
        \begin{equation}
            \dt \sint{T} = -\int_{\Gamma}\vec{\mathbf{F}}\dotn \,ds + W\big\vert_{\partial\mathscr{S}}. \label{eqn:conlawDef}
        \end{equation}
        \Cref{eqn:conlawDef} expresses the fact that the rate of changes in total density equals the sum of the fluxes over the fixed boundary \(\Gamma\) and the boundary of the free surface \(\partial\mathscr{S}\). 
        
        In \cite{olver}, the conserved densities are first found via symmetry methods.  Then, these densities (shown in \Cref{tab:summary}) are validated by differentiating the following contour integrals with respect to time: 
        \begin{subequations}
            \begin{align}
                &I_1 = -\oint_{\partial\mathscr{D}} \phi \,dz = \dint{-\phi_x}\label{eqn:I1}\\
                &I_2 = \oint_{\partial\mathscr{D}} \frac{1}{2}\phi\grad\phi\dotn\,ds + \frac{1}{2}gz^2\,dx = \dint{\frac{1}{2}\left\vert \grad\phi\right\vert^2 + gz}\label{eqn:I2}\\
                &I_3 =  \oint_{\partial\mathscr{D}} z\,dx = \iint_{\mathscr{D}}\,dx\,dz\label{eqn:I3}\\
                &I_4 = \oint_{\partial\mathscr{D}}\phi\,dx=\dint{\phi_z}\label{eqn:I4}\\
                &I_5 = \oint_{\partial\mathscr{D}} xz\,dx = \dint{x}\label{eqn:I5}\\
                &I_6 = \oint_{\partial\mathscr{D}} \frac{1}{2}z^2\,dx = \dint{z}\label{eqn:I6}\\
                &I_7 = \oint_{\partial\mathscr{D}} \phi\left(z\,dx - x\,dz\right) = \dint{x\phi_x + z\phi_z + 2\phi}\label{eqn:I7}\\
                &I_8 = \oint_{\partial\mathscr{D}} \phi\left(x\,dx + z\,dz\right) = \dint{x\phi_z - z\phi_x}\label{eqn:I8}.
            \end{align}
        \end{subequations}
        In the following section, we demonstrate how the above eight contour integrals directly arise from the weak formulation presented in \cite{NonlocalNonlocal}.  

%% file: sections/equations.tex
        There are many different reformulation of \eqref{eqn:laplace1d}-\eqref{eqn:dynamic1d} that yield the equations of motion in terms of free-surface, and the trace of the velocity potential along the free-surface.   For a summary, see \cite{wilkening2015comparison} for the formulations.

        One such reformulation serves as a starting point for the work presented here.  The formulation given by Ablowitz, Fokas, \& Musslimani \cite{ablowitz2006new} yields a coupled system of differential and integro-differential equations for \(\eta(x,t)\) and \(q(x,t) = \phi(x,\eta(x,t),t)\).  In their work, the authors begin by considering a harmonic test function \(\varphi(x,z)\) and the velocity potential \(\phi\).  Via a divergence theorem argument, they recast the kinematic boundary condition \eqref{eqn:kinematic1d} as 
        \begin{equation}
            \cint{\varphi_z\left(\grad\phi\dotn\right) - \varphi_x\left(\grad\phi\dott\right)} = 0, \label{eqn:GreenPhiAFM}
        \end{equation}
        where \(\vn\) is the outward pointing normal, \(\vt\) is the tangent vector, and both \(\eta\) and \(\vert\grad\phi\vert\) decay sufficiently fast as \(|x|\to\infty\) such that the following integrals make sense.  Choosing \(\varphi = e^{-ikx}\sinh(k(z+h))\) eliminates any contribution from the bottom boundary and yields the system written in terms of the surface variables \(\eta(x,t)\) and \(q(x,t) = \phi(x,\eta(x,t),t)\) are given by
        \begin{eqnarray}
            &\displaystyle \wint{e^{-ikx}\left(\eta_t\cosh(k(\eta+h)) + iq_x\sinh(k(\eta+h))\right)} = 0,\label{eqn:AFMNonlocal}\\
            &\displaystyle q_t + \frac{1}{2}q_x^2 + g\eta - \frac{1}{2}\frac{(\eta_t + q_x\eta_x)^2}{1 + \eta_x^2} = 0.\label{eqn:AFMLocal}
        \end{eqnarray}
        One can think of \eqref{eqn:AFMNonlocal} as an implicit relationship defining the Dirichlet-to-Neumann map at the free surface.  Indeed, Haut \& Ablowitz \cite{haut_ablowitz_2009} show that \eqref{eqn:AFMNonlocal} can be used to generate the same expansion of the DNO as found by Craig \& Sulem \cite{craig1993numerical}. The authors also establish the equivalence of the two formulations with that of Zakharov \cite{zakharov1968stability} and Craig \& Sulem \cite{craig1993numerical}; see \cite{haut_ablowitz_2009} for details.

        Recently, an extended version of the local/nonlocal formulation was introduced in \cite{NonlocalNonlocal} by considering a harmonic test function \(\varphi\) and by noting that both \(\phi\) and \(\phi_t\) are harmonic within the fluid bulk.  Via the appropriate integral theorems, and by assuming sufficient decay in \(\eta\) and \(\phi\), we begin by considering the integral relations
            \begin{eqnarray}
                \cint{\varphi_z\left(\grad\phi\dotn\right) - \phi\left(\grad\varphi_z\dotn\right)} = 0, \label{eqn:boundA}\\
                \cint{\varphi_z\left(\grad\phi_t\dotn\right) - \phi_t\left(\grad\varphi_z\dotn\right)} = 0 \label{eqn:boundB}.
            \end{eqnarray}
        As discussed in \cite{NonlocalNonlocal}, the resulting equations in terms of boundary variables are given by \begin{eqnarray}
            &\displaystyle \sint{\varphi_z\eta_t - q\left(\grad\varphi_z\dotn\right)} = \bint{-Q\,\varphi_{zz}},\label{eqn:nlnlA}\\
            &\displaystyle \sint{q_x\eta_t\varphi_{xz} - \left(q_t + g\eta\right)\left(\varphi_{zz} + \eta_x\varphi_{xz}\right)} = \bint{\frac{1}{2}Q_x^2\,\varphi_{zz}}\label{eqn:nlnlB},
        \end{eqnarray}
        where we have introduced \(Q(x,t) = \phi(x,-h,t)\) to represent the trace of the velocity potential along the bottom of the fluid. Equations \eqref{eqn:nlnlA} and \eqref{eqn:nlnlB} are valid for any harmonic function \(\varphi(x,z)\) along with suitable decay  properties for \(\eta\) and \(\vert\grad\phi\vert\) such that the above integrals make sense.  Here, we have also introduced the notation \[\int_{\mathscr{S}}f(x,z,t)\,dx = \int_\mathbb{R} f(x,\eta,t)\,dx \quad \text{and}\quad \int_{\mathscr{B}}f(x,z,t)\,dx = \int_{\mathbb{R}}f(x,-h,t)\,dx\] to denote the integral over the whole line where \(z = \eta(x,t)\) and \(z = -h\) respectively.
        
        \begin{remark}
            It is worth noting that the reformulation \eqref{eqn:nlnlA}-\eqref{eqn:nlnlB} is not completely written in terms of surface variables as the velocity potential along the bottom, \(Q(x,t)\), still remains a part of the formulation.  However, as noted in \cite{NonlocalNonlocal}, choosing a family of test functions in the form \(\varphi(x,z) = e^{-ikx}\sinh(k(z+h))\) eliminate any contribution from the bottom terms so that \eqref{eqn:nlnlA} - \eqref{eqn:nlnlB} represent a solvable system in terms of the surface variables \(\eta(x,t)\) and \(q(x,t)\).  However, in order to derive the conservation laws, we examine a variety of different harmonic functions \(\varphi\).  
        \end{remark}
        

%% file: sections/conservation_laws_1D_surface.tex

    To find conservations laws directly from the equations of motion, we begin by rewriting \eqref{eqn:nlnlA} and \eqref{eqn:nlnlB} as 
    \begin{eqnarray}
        &&\dt\sint{\partial_z^{-1}\varphi} = \oint_{\partial\mathscr{D}}\phi\left(\sigma_3\grad\varphi\right)\dotn\,ds + 2\odint{\phi\varphi_{z}\,dx}\label{eqn:nlnlAA}\\ 
        \nonumber\\
        &&\dt\oint_{\partial\mathscr{D}}{\left(\phi\left(\sigma_3\grad\varphi\right)\dotn\right)\,ds}=-\oint_{\partial\mathscr{D}}{\left[(p/\rho + gz)\left(\sigma_3\grad\varphi\right)\dotn +2\phi\varphi_{zz}(\grad\phi\dotn)\right]\,ds},\label{eqn:nlnlBB}
    \end{eqnarray}
    where \(\partial_z^{-1}\varphi\) is the anti-derivative of \(\varphi\) with respect to \(z\) such that \(\partial_z^{-1}\varphi(x,z)\bigg\vert_{z = 0} = 0\) for all \(x\) and \(\sigma_3\) is the third Pauli matrix so that 
    \[
        \sigma_3\grad\varphi = \begin{bmatrix}1 & 0\\0&-1\end{bmatrix}\begin{bmatrix}\varphi_{x}\\\varphi_{z} \end{bmatrix}= \begin{bmatrix}\varphi_{x}\\-\varphi_{z}\end{bmatrix}.
    \] 
    In \eqref{eqn:nlnlAA} and \eqref{eqn:nlnlBB}, we have assumed that \(\grad\phi\dotn = 0\) on \(\Gamma\) as described in \Cref{sec:background}.
    
    With this simplification, we now proceed to compute the eight conservation laws as computed by \cite{olver} by choosing a sequence of harmonic polynomials in the form \[\varphi_n = \frac{1}{n!}(x + iz)^n.\]  As noted in \cite{olver}, there is a recursive nature where \emph{higher-order} conserved densities are found in terms of lower-order densities.  With this in mind, we define the quantities \(A\) and \(B\) as 
    \begin{equation}
        A = \sint{\partial_z^{-1}\varphi}, \qquad B = \odint{\phi(\sigma_3\grad\varphi)\dotn}\,ds.\label{eqn:ABEqns}
    \end{equation}
    With these definitions, Equations \eqref{eqn:nlnlAA} and \eqref{eqn:nlnlBB} can be represented in the following compact form that is conducive to back-substituting lower-order conserved densities.
    \begin{eqnarray}
        \frac{dB}{dt} &=& -\odint{(p/\rho + gz)(\sigma_3\grad\varphi)\dotn}\,ds - 2\sint{\phi\,\eta_t\varphi_{zz}},\label{eqn:nlnlBBB}\\
        \frac{dA}{dt} &=& B + 2\odint{\phi\varphi_{z}}\,dx. \label{eqn:nlnlAAA}
    \end{eqnarray}
    
    \Cref{tab:summary} shows a summary of the conserved densities at the free surface along with the corresponding test functions.  The boundary fluxes are given in the details presented below.
    \subsection{\(\varphi_0 = 1\)}
        Equations \eqref{eqn:nlnlBBB} and \eqref{eqn:nlnlAAA} serve as our foundation for deriving the conservation laws. As stated earlier, we assume that \(\phi\) and \(\eta\) decay rapidly enough as \(\vert x \vert \to \infty\) such that the corresponding integrals make sense.  Substituting \(\varphi_0 = 1\) into the expressions for \(A\) and \(B\) yields 
        \[A_0 = \sint{\eta} \qquad B_0 = 0.\]
        Using \eqref{eqn:nlnlBBB} yields \(0 = 0\), while substitution into \eqref{eqn:nlnlAAA} yields 
        \begin{equation}
            \dt A_0 = \dt\sint{\eta} = 0 \qquad \Rightarrow\qquad \dt\wint{\eta} = 0. \qquad \label{eqn:T3}\tag{\(T_3\)}
        \end{equation}
        where we have used the same numbering convention for the conserved densities as in \cite{olver}.  
    \subsection{\(\varphi_1 = x + iz\)}
        Substituting \(\varphi = x + iz\) into the expressions for \(A\) and \(B\) given by \eqref{eqn:ABEqns} we find 
        \[A_1 = \sint{\left(x\eta + \frac{i}{2}\eta^2\right)}\,dx \qquad B_1 = \odint{\phi\begin{bmatrix}1\\-i\end{bmatrix}\dotn}\,ds.\]
        Using \eqref{eqn:nlnlBBB} we find 
        \begin{eqnarray*}
            \dt B_1 
            &=& igA_0-i\bint{p/\rho} + g\odint{z}\,dz\\
            &=& i\left(gA_0 - \bint{p/\rho}\right)\\
            &=& i\dt\left(tgA_0\right)- i\bint{p/\rho} - igt\dt\left(A_0\right)
        \end{eqnarray*}
        
        As previously stated, we assume that \(\eta\to 0\) rapidly so that \(\eta^2\) vanishes in the limit as \(|x|\to\infty\).  Using the fact that \(\dt A_0 = 0\) from \eqref{eqn:T3}, the real and imaginary parts yield the following conservation laws with the appropriate boundary-flux terms denoted: 
        \begin{equation}
            \dt\odint{-\phi}\,dz = 0,\label{eqn:T1}\tag{\(T_1\)}
        \end{equation}
        and 
        \begin{equation}
            \dt\odint{-\phi - tgz}\,dx = -\bint{p/\rho}.\label{eqn:T4}\tag{\(T_4\)}
        \end{equation}
        
        It is worth noting that in deriving \eqref{eqn:T1} and \eqref{eqn:T4}, the previously determined conservation law \eqref{eqn:T3} was needed.  This will be a recurring theme as we move up the hierarchy; conservation laws determined from lower-order harmonic polynomials will be needed to move to higher-order harmonic polynomials. 

        Repeating the same process with \(A_1\), and using the integral identities generated in \eqref{eqn:nlnl:idB}, upon separating real and imaginary parts, we find  
        \begin{equation}
           \displaystyle \dt\odint{ xz\,dx + t\phi\,dz} = 0,\label{eqn:T5}\tag{\(T_5\)}
        \end{equation}
        and 
        \begin{equation}
            \dt\sint{-\frac{1}{2}\eta^2} +  \dt\odint{ t \phi\,dx + \frac{1}{2}gt^2 z\,dx} = 0.\label{eqn:T6}\tag{\(T_6\)}
        \end{equation}
    \subsection{\(\varphi_2 = \frac{1}{2}\left(x + iz\right)^2\)}
        To continue in the hierarchy, we continue the process by beginning with \(\varphi_2 = \frac{1}{2}(x + iz)^2\) substituted into \eqref{eqn:nlnlBBB}.  Separating real and imaginary parts, we find
        \begin{equation}
            \dt\odint{\phi(z\,dx - x\,dz) -\left(4t\mathscr{H} -\frac{7}{2}gz^2 - \frac{7}{6}g^2t^3 z\right)\,dx} = \bint{hp/\rho}, \label{eqn:T7}\tag{\(T_7\)}
        \end{equation}
        \begin{equation}
            \dt\odint{\phi(x\,dx + z\,dz) + gtxz\,dx + \frac{1}{2}gt^2\phi\,dz} = -\bint{x p/\rho },\label{eqn:T8}\tag{\(T_8\)}
        \end{equation}
        where we have used the fact that \[\dt \mathcal{H} = \dt\sint{\mathscr{H}} = \dt\sint{\underbrace{\frac{1}{2}\left(\phi\grad\phi\dotn + gz^2\right)}_{\mathscr{H}}} = 0.\]  
        \begin{remark}
            Knowing the Hamiltonian at the outset does simplify the calculations.  However, a priori knowledge of \(\mathscr{H}\) is not necessary.  Instead, the Hamiltonian arises naturally when one considers the integral identities presented in \Cref{app:simplification1d}.
        \end{remark}
    \subsection{Summary of Conservation Laws}
        At this point, we have computed the eight conservation laws presented in \cite{olver} for the case with no surface tension.  The corresponding densities at the free surface are given in \Cref{tab:summary}.  The corresponding boundary flux terms can be deduced from Equations \eqref{eqn:T1} - \eqref{eqn:T8}.
        \begin{table}[H]
            \hrule ~\\
            \caption{Table of Conserved Densities from \cite{olver} (Thm. 6.2) along with a summary of how they were found via the weak formulation.}\label{tab:summary}
            
            \vspace*{.05in}
            \renewcommand{\arraystretch}{1.8}
            \hrule
            \small{
            \begin{tabular}{l|l}
                \(\displaystyle T_{1,s} = -\eta_x q\)              & Imaginary part of \eqref{eqn:nlnlBB} with \(\varphi = \frac{1}{2}(x + iz)^2\)\\
                \(\displaystyle T_{2,s} = \frac{1}{2}q\eta_t + \frac{1}{2}g\eta^2\) & Embedded in \eqref{eqn:nlnlBB} with \(\varphi = \frac{1}{6}(x+iz)^3\)\\
                \(\displaystyle T_{3,s} = \eta\) & Imaginary part of \eqref{eqn:nlnlAA} with \(\varphi = x +iz\)\\
                \(\displaystyle T_{4,s} = q + gt\eta\) &Real part of \eqref{eqn:nlnlBB} with \(\varphi =\frac{1}{2} (x + iz)^2\)\\ 
                \(\displaystyle T_{5,s} = x\eta + t\eta_x q\)  &Imaginary part of \eqref{eqn:nlnlAA} with \(\varphi = \frac{1}{2}(x + iz)^2\)\\
                \(\displaystyle T_{6,s} = \frac{1}{2}\eta^2 - tq - \frac{1}{2}gt^2\eta\)  &Real part of \eqref{eqn:nlnlAA} with \(\varphi = \frac{1}{2}(x+iz)^2\)\\
                \(\displaystyle T_{7,s} = q(\eta - x\eta_x) - 4tT_2 + \frac{7}{2}gt\eta^2 - \frac{7}{2}gt^2q - \frac{7}{6}g^2t^3\eta\)   &Imaginary part of \eqref{eqn:nlnlBB} with \(\varphi = \frac{1}{6}(x + iz)^3\)\\
                \(\displaystyle T_{8,s} = (x + \eta\eta_x)q + gtx\eta +\frac{1}{2}t^2g \eta_x q\) &Real part of  \eqref{eqn:nlnlBB} with \(\varphi = \frac{1}{6}(x + iz)^3\)\\
            \end{tabular}}
            \hrule
        \end{table}

        As in \cite{olver,olver1983conservation}, deriving \eqref{eqn:T7} requires additional manipulations.  While we do not need to take specific linear combinations of separate one-forms, we do instead need to find and extract the Hamiltonian.  This is aided by the observation of the term \(2\eta_t\,\phi\,\varphi_{zzz}(x,\eta)\) as a \emph{Legendre-type transform} lurking within Equation \eqref{eqn:nlnlBBB}.  At this point, we have derived the 8 conservation laws for a one-dimensional free-surface on the whole-line (see \Cref{tab:summary} for details) without surface-tension.  The computed values for \(T_1\) - \(T_8\) are in agreement with those presented in \cite{olver,olver1983conservation}.  There is nothing that necessarily stops us from continuing the process.  For example, we have yet to use the cubic test function in \eqref{eqn:nlnlAA}.  However, without introducing nonlocal operators of some variety, it appears to be impossible to continue the aforementioned process as expected via the results of \cite{olver1983conservation}.  
    
        It is a rather remarkable fact that by considering the harmonic polynomials of the form \[\varphi_n = \frac{1}{n!}(x + iz)^n,\] we are able to immediately recover the various contour integrals used to prove Theorem 6.1 in \S6 of \cite{olver}.  Specifically, \(I^1, I^3, I^4, \ldots, I^8\) are found via the left-hand sides of Equations \eqref{eqn:nlnlAA} - \eqref{eqn:nlnlBB}, while the closed differentials \(G^1, G^3, G^4, \ldots, G^8\) are found directly via \eqref{eqn:nlnl:idA} and \eqref{eqn:nlnl:idB} as described in \Cref{app:simplification1d}.  

    \subsection{Including Surface Tension}
        When surface tension is included, the dynamic boundary condition \eqref{eqn:dynamic1d} is replaced by 
        \begin{equation}
            \phi_t + \frac{1}{2}\vert\grad\phi\vert^2 + g\eta =\frac{\sigma\eta_{xx}}{\left(1 + \eta_x^2\right)^{3/2}}\label{eqn:dynamic1DST},
        \end{equation}
        where \(\sigma\) represents the coefficient of surface-tension.  The resulting formulations still hold with the slight modification as presented below:
        \begin{eqnarray}
            &&\dt\odint{\partial_z^{-1}\varphi}\,dx = \oint_{\partial\mathscr{D}}\phi\grad\phi\dotn\,ds\label{eqn:nlnlAAST}\\ \nonumber\\
            &&\dt\oint_{\partial\mathscr{D}}{\left(\phi\left(\sigma_3\grad\varphi\right)\dotn\right)\,ds}=-\oint_{\partial\mathscr{D}}{\left[(p/\rho + gz + S)\left(\sigma_3\grad\varphi\right)\dotn +2\phi\varphi_{zz}(\grad\phi\dotn)\right]\,ds},\label{eqn:nlnlBBST}
        \end{eqnarray}
        where \[S = -\dx\frac{\sigma\eta_x}{\sqrt{1+\eta_x^2}}\] represent the surface tension force acting only at the free surface \(\mathscr{S}\).  
        
        We find the same results as described in \Cref{tab:summary} with the exception of two changes. First, the resulting Hamiltonian or \(T_2\) is modified to include the additional term 
        \[T_2 = \frac{1}{2}q\eta_t + \frac{g}{2}\eta^2 + \sigma\left(\sqrt{1 + \eta_x^2} - 1\right).\]  Second, when attempting to recover \eqref{eqn:T7}, one encounters an additional term of the form \[\dt\sint{\sigma\sqrt{1+\eta_x^2}}\neq 0.\]  This results in the loss of \(T_7\) as a conserved density in the presence of surface tension and is to be expected (see \cite{olver, olver1983conservation}).  Physically, it makes sense as to why \eqref{eqn:T7} is no longer conserved in the presence of surface tension as it does not make sense for the arc-length of the interface to remain constant for all time.
        

%% file: sections/conservation_laws_linear_shear.tex
\begin{figure}[htp]
  \centering
  \begin{tikzpicture}
    \begin{axis}[
      hide axis,
      xmin        =   -13,
      xmax        =   13,
      ymin        =   -1.15,
      ymax        =   .8,
      domain      =   -10:10,
      width       =   \textwidth,
      height      =   1.5in,
      samples     =   300,
      black,
      thick
  ]
    \addplot[name path = wave] {.35*2^(-.2*x^2)} node [pos=1, right] {\(z = \eta(x,t)\)}; 
    \addplot[name path = bottom]{-1} node [pos = 1, right] {\(z = -h\)};
    \addplot[name path = floor] {-1.15}; 

    \addplot[fill = NavyBlue!30] fill between[of=wave and bottom];
    \addplot[fill=Brown!60] fill between[of=floor and bottom];

      \addplot[color=black,thick,->] coordinates {
        (-6,-.8)
        (-8.5,-.8)
        };
      \addplot[color=black,thick,->] coordinates {
        (-6,-.6)
        (-8,-.6)
        };
      \addplot[color=black,thick,->] coordinates {
        (-6,-.4)
        (-7.5,-.4)
        };

      \addplot[color=black,thick,->] coordinates {
        (-6,-.2)
        (-7,-.2)
        };
      \addplot[color=black,thick,dashed] coordinates {
        (-6,0)
        (-6,-1)
        };
      \addplot[color=black,thick,dashed] coordinates {
        (-6.5,0)
        (-9,-1)
        };

    \end{axis}
  \end{tikzpicture}
  \caption{Fluid Domain for Constant Vorticity}\label{fig:fluidDomainCV}
\end{figure}
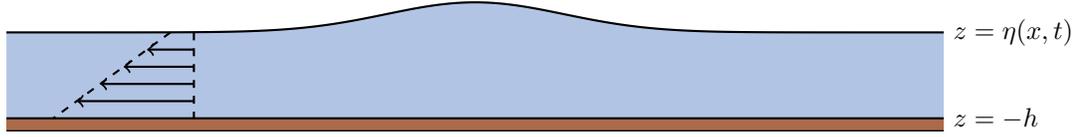

As discussed in \cite{ashton2011non}, the formulation given by \cite{ablowitz2006new} can be extended to include constant vorticity.  Here, we extend the nonlocal-nonlocal formulation of \cite{NonlocalNonlocal} to include constant vorticity and examine the conservation laws in the presence of linear shear.  

We begin by considering the equations of motion expressed in terms of the fluid velocities given by the following bulk equations
\begin{align}
  &\grad\cdot\vec{\mathbf{u}} = 0, \quad \grad \times\vec{\mathbf{u}} = \omega, &&(x,z)\in \mathscr{D},\label{eqn:conMassConVort}\\
  &u_t + uu_x + v u_z + \frac{p_x}{\rho} = 0, &&(x,z)\in \mathscr{D}\label{eqn:uB}\\
  &v_t + uv_x + v v_z + g + \frac{p_z}{\rho}, &&(x,z)\in \mathscr{D}\label{eqn:vB}
\end{align}
with corresponding boundary conditions: 
\begin{align}
  &v = 0, && z= -h\label{eqn:kinematicBottomCV}\\
  &\eta_t = v - \eta_x u && z= \eta,\label{eqn:kinematicCV}\\
  &p = 0, && z= \eta, \label{eqn:pressureBCCV}
\end{align} 
where \(\omega\) is the vorticity and remains constant throughout the fluid,  and \(\vec{\mathbf{u}} = \begin{bmatrix} u & v\end{bmatrix}^{\textsf{T}}\) with \(u\) and \(v\) representing the horizontal and vertical velocities of the fluid respectively.  Following \cite{wahlen2007hamiltonian,ashton2011non}, we can introduce the pseudo-potential function \(\phi\) such that 
\begin{equation}
  \grad\phi = \begin{bmatrix}u + \omega z\\v\end{bmatrix}\label{eqn:pseudoPot}.
\end{equation}

Again, we assume that \(\grad\phi\) and \(\eta\) decay sufficiently fast enough as \(\vert x \vert \to \infty\).
Following the work of \cite{gavrilyuk2015kinematic}, we note that the free-surface and bottom both correspond to streamlines.  This allows us to easily re-write a total \(x\) derivative of the momentum equations along the free-surface and bottom of the fluid.  Substituting \eqref{eqn:pseudoPot} into \eqref{eqn:conMassConVort} - \eqref{eqn:pressureBCCV} yields the following equivalent problem in terms of the pseudo-potential \(\phi\):
\begin{align}
  &\phi_{xx} + \phi_{zz} = 0, \qquad &&(x,z)\in\mathscr{D},\\
  &\phi_z = 0, &&z = -h,\\
  &\phi_t + \frac{1}{2}\phi_x^2 + \omega h\phi_x - gh + \frac{p}{\rho} = M_\mathscr{B}(t), && z = -h,\\
  &\eta_t = \phi_z - \eta_x\phi_x + \omega\eta\eta_x, && z = \eta(x,t),\\
  &\phi_t + \frac{1}{2}\phi_x^2 + \frac{1}{2}\phi_z^2 - \omega\eta\phi_x + \frac{\omega^2}{2}\eta^2 - \omega\partial_x^{-1}\eta_t + g\eta = M_\mathscr{S}(t),&&z = \eta(x,t),
\end{align}
where \(M_\mathscr{B}\) and \(M_{\mathscr{S}}\) represent an arbitrary function of integration that depends solely on \(t\).  Contrary to the irrotational case, this function need not be the same along each streamline, and thus, we prescribe separate functions for each streamline.  Without loss of generality, we take \(M_\mathscr{S}(t) = 0\) and leave \(M_\mathscr{B}(t)\) arbitrary.  Now that the problem is formulated in terms of a harmonic function \(\phi\), we can proceed as before. 

\subsection{Derivation of the Nonlocal/Nonlocal Formulation \& Conservation Laws}
  Following the same process as outlined in \Cref{sec:conlaw_irrotational} using the modified boundary value problem described above yields the corresponding nonlocal/nonlocal formulation in a form conducive for computing conservation laws:
  \begin{equation}
    \displaystyle\sint{\dt\varphi} = \sint{\omega\eta\eta_x\varphi_z} + \oint_{\partial\mathscr{D}}{\phi\,\grad\varphi_z\dotn\,ds},\label{eqn:vortnlnlAA}
  \end{equation}
  along with 
  \begin{eqnarray}
    \dt B =   \sint{-g\eta\left(\sigma_3\grad\varphi_z\right)\dotn-\omega\eta\left(q_x - \frac{\omega\eta}{2}\right)\varphi_{zz} -\omega\eta\left(\eta_t - \omega\eta\eta_x\right)\varphi_{xz}- 2\left(\phi - \omega\partial_x^{-1}\eta\right)\eta_t\varphi_{zzz}}\nonumber\\
    +\bint{\frac{1}{2}Q_x^2\varphi_{zz}}\label{eqn:vortnlnlBB},
  \end{eqnarray}
  where \(B\) is given by \[B = \sint{\phi\left(\sigma_3\grad\varphi_z\right)\dotn}.\]
  As before, we define \(q\) and \(Q\) to be the trace of the pseudo-potential along the free-surface and bottom respectively.  That is, \[q(x,t)=\phi(x,\eta,t) \quad \text{and}\quad Q(x,t) = \phi(x,-h,t).\] This choice of variables are convenient to use and have specific advantages.  However, one disadvantage is that \(q\) and \(\eta\) are no longer the canonical Hamiltonian variables for this problem.  As shown in \cite{wahlen2007hamiltonian}, the appropriate canonical Hamiltonian variables are given by \((\zeta, \eta)\) where 
  \begin{equation}
    \zeta = q - \frac{\omega}{2}\partial_x^{-1}\eta.\label{eqn:CVcanVar}
  \end{equation}

  By taking the limit as \(\omega\to 0\) and performing the appropriate simplifications, Equations \eqref{eqn:vortnlnlAA} and \eqref{eqn:vortnlnlBB} reduce to \eqref{eqn:nlnlAA} and \eqref{eqn:nlnlBB} respectively. 

  For brevity, we omit the computations for the conservation laws with constant vorticity as they follow the same process in \Cref{sec:conlaw_irrotational} and do not shed any additional light on the problem.  Instead, we present the results in \Cref{tab:summaryConstantVorticity} in terms of the canonical Hamiltonian variables \((\zeta,\eta)\) (with the exception of the Hamiltonian) as defined in \eqref{eqn:CVcanVar}.

  \begin{table}
    \hrule ~\\
    \vspace*{.05in}
    \caption{Table of Conserved Densities for constant vorticity}\label{tab:summaryConstantVorticity}
    \renewcommand{\arraystretch}{1.8}
    \hrule
    \small{
    \begin{tabular}{l|l}
      \(\displaystyle T_1 = -\eta_x \zeta\)              & Imaginary part of \eqref{eqn:vortnlnlBB} with \(\varphi = \frac{1}{2}(x + iz)^2\)\\
      \(\displaystyle T_2 = \frac{1}{2}\phi\eta_t +\omega\eta\eta_x\phi + \frac{\omega^2}{6}\eta^3 + \frac{g}{2}\eta^2\) & Embedded in \eqref{eqn:vortnlnlBB} with \(\varphi = \frac{1}{6}(x+iz)^3\)\\
      \(\displaystyle T_3 = \eta\) & Imaginary part of \eqref{eqn:vortnlnlAA} with \(\varphi = x +iz\)\\
      \(\displaystyle T_4 = \zeta + \frac{\omega}{2}x\eta + gt\eta + \omega t \zeta\eta_x\) &Real part of \eqref{eqn:vortnlnlBB} with \(\varphi =\frac{1}{2} (x + iz)^2\)\\ 
      \(\displaystyle T_5 = x\eta + t\zeta\eta_x \)  &Imaginary part of \eqref{eqn:vortnlnlAA} with \(\varphi = \frac{1}{2}(x + iz)^2\)\\
      \(\displaystyle T_6 = \frac{1}{2}\eta^2 - t\zeta  - \frac{gt^2}{2}\eta+ \frac{\omega x t}{2}\eta\)  &Real part of \eqref{eqn:vortnlnlAA} with \(\varphi = \frac{1}{2}(x+iz)^2\)\\
      \(\displaystyle T_7 =\) \emph{\bf undetermined}   &Imaginary part of \eqref{eqn:vortnlnlBB} with \(\varphi = \frac{1}{6}(x + iz)^3\)\\
      \(\displaystyle T_8 = \zeta(x + \eta\eta_x)  + xgt\eta + \frac{gt^2}{2}\eta_x\zeta - \frac{\omega x^2\eta}{4} + \frac{\omega}{12}\eta^3 \) &Real part of  \eqref{eqn:vortnlnlBB} with \(\varphi = \frac{1}{6}(x + iz)^3\)\\[1em]
    \end{tabular}}
    \hrule
  \end{table}
  It is worth noting that Longuet-Higgins \cite{longuetIntegrals} extended seven of Benjamin \& Olver's integral identities to a fully rotational setting.  In \S6 of \cite{longuetIntegrals}, the following integrals are introduced:
  \begin{subequations}
      \begin{align}
          &I_1^*   = \dint{u}\label{eqn:rot:I1}\\
          &I_2^*   = \dint{\frac{1}{2}\left(u^2 + v^2\right) + gz}\label{eqn:rot:I2}\\
          &I_3^*   = \iint_{\mathscr{D}}\,dx\,dz\label{eqn:rot:I3}\\
          &I_4^*   = \dint{v}\label{eqn:rot:I4}\\
          &I_5^*   = \dint{x}\label{eqn:rot:I5}\\
          &I_6^*   = \dint{z}\label{eqn:rot:I6}\\
          &I_7^* = \dint{v_x - u_z}\label{eqn:rot:I7}\\
          &I_8^*   = \dint{xv - zu}\label{eqn:rot:I8}.
      \end{align}
  \end{subequations}
  In the above, \(I_7^*\) is a replacement for \(I_7\) and represents the circulation.  However, we can easily write a more general version such that for any differentiable function \(f(\omega)\), the quantity \[\dt\iint_\mathscr{D}f(v_z - u_x)\,dz\,dx = 0\] would also yield a conservation law. 
  
  The author presented the above integrals and proceeded to show that these integrals eventually yielded conserved quantities provided that the pressure vanish on the boundary.  However, these integral identities were presented and confirmed rather than derived.  
  
  While \(I_1^*, I_2^*, \ldots, I_8^*\) are conserved, it remains unclear if the variables can be returned from the bulk to the boundary.  This can be compared with the results found for constant vorticity.  While we did not present any new identities, we have shown (1) that the integral identities can be derived directly from the equations of motion for constant vorticity, and (2) how the conserved densities can be expressed in terms of surface variables the sum of the fluxes over the fixed boundary \(\Gamma\).  

%% file: sections/acknolwledgements.tex
The authors would like to Vishal Vasan for fruitful discussions related to this work and Miles Wheeler for pointing us towards \cite{longuetIntegrals}.  

Portions of this work were conducted while KO was located at ICERM at Brown University, University of Vienna, and Yale-NUS.  Both KO and SCY gratefully acknowledge the support of the National Science Foundation under Grant Number DMS-1715082. Any opinions, findings, and conclusions or recommendations expressed in this material are those of the author and do not necessarily reflect the views of the National Science Foundation.

%% file: sections/appendix.tex
\appendix
    \section{Simplification Formulae for Irrotational Fluids}\label{app:simplification1d}
    In the derivation of \eqref{eqn:nlnlB} (and thus, \eqref{eqn:nlnlBB}), for any differentiable function \(f(\phi_x + i\phi_z) = U(\phi_x,\phi_z) + iV(\phi_x,\phi_z)\) we encountered a hierarchy of additional integral relationships given by
    \begin{equation}
        \odint{U(\phi_x,\phi_z)\grad\varphi_z\dotn + V(\phi_x,\phi_z)\grad\varphi_x\dotn} = 0,\label{eqn:app:nlnlACA}
    \end{equation}
    or in an alternative form, 
    \begin{equation}
        \odint{ \varphi_{z}\begin{bmatrix}U\\V\end{bmatrix}\dotn - \varphi_{x}\begin{bmatrix}-V\\~U\end{bmatrix}\dotn +}\,ds = 0.\label{eqn:app:genGreenPhiAFM}
    \end{equation}
    Equations \eqref{eqn:app:nlnlACA} and \eqref{eqn:app:genGreenPhiAFM} is simply a generalization of \eqref{eqn:GreenPhiAFM} and are useful when computing the fluxes.  In fact, by choosing 
    \[
        f_1(\phi_x + i\phi_z) = \phi_x + i\phi_z, \quad \text{and}\quad f_2(\phi_x + i\phi_z) = \frac{1}{2}(\phi_x + i\phi_z)^2,
    \]
    we find the following relationships that are useful for eliminating boundary terms:
    \begin{eqnarray}
        &\displaystyle\odint{\phi_x\grad\varphi_z\dotn + \phi_z\grad\varphi_x\dotn} = 0,\label{eqn:nlnl:idA}\\
        &\displaystyle\odint{\frac{1}{2}\left(\phi_x^2 - \phi_z^2\right)\grad\varphi_z\dotn + \phi_x\phi_z\grad\varphi_x\dotn } = 0.\label{eqn:nlnl:idB}
    \end{eqnarray}
    These are precisely the same vanishing contour integrals labeled \(G_j\) in \S6 of \cite{olver} for the functions \(\varphi = \frac{1}{n!}(x + iz)^n\) for \(n = 1, 2\).  
    
    It is worth noting that we can continue this hierarchy, though the usefulness of the following have yet to be explored.  For example, if \(f_3\) is chosen so that \[f_3(\phi_x + i\phi_z) = \frac{1}{3}(\phi_x + i\phi_z)^3,\] then 
    \begin{equation}
        \displaystyle\odint{\left(\frac{1}{3}\phi_x^3 - \phi_x\phi_z^2\right)\grad\varphi_z\dotn + \left(\phi_x^2\phi_z - \frac{1}{3}\phi_z^3\right)\grad\varphi_x\dotn} = 0.\label{eqn:nlnlDD}
    \end{equation}
    If \(\varphi = xz\), then we find
    \begin{equation}
        \odint{\left(\frac{1}{3}\phi_x^3 - \phi_x\phi_z^2\right)(-1)\,dz + \left(\phi_x^2\phi_z - \frac{1}{3}\phi_z^3\right)\,dx} = 0.\label{eqn:app:third}
    \end{equation}
    Of course, without additional boundary data related to the problem of interest, \eqref{eqn:nlnl:idA} -\eqref{eqn:nlnlDD} are just specific examples of Green's identities.  However, once the boundary conditions given by  \eqref{eqn:kinematicBottom1d} - \eqref{eqn:dynamic1d} are incorporated, then the resulting integral relationships are useful in constructing the conservation laws.  